\newtheorem{lemma}{Lemma}

\newtheorem{theorem}[lemma]{Theorem}

\newtheorem{rem}[lemma]{Remark}
\newtheorem{frem}[lemma]{Final remark}

\newcommand{\kla}{\left ( }
\newcommand{\nach}{\rightarrow}
\newcommand{\mer}{\right ) }
\newcommand{\qed}{\hspace*{\fill}$\Box$\hz\pagebreak[1]}

\newcommand{\for}{\begin{eqnarray*}}
\newcommand{\mel}{\end{eqnarray*}}

\newcommand{\mitt}{\left | { \atop } \right.}
\newcommand{\kl}{\pl \le \pl}

\newcommand{\kll}{\p \le \p}
\newcommand{\gll}{\p \ge \p}

\newcommand{\nz}{{\rm  I\! N}}
\newcommand{\nen}{n \in \nz}

\newcommand{\rz}{{\rm  I\! R}}

\newcommand{\p}{\hspace{.05cm}}
\newcommand{\pl}{\hspace{.1cm}}
\newcommand{\pll}{\hspace{.3cm}}
\newcommand{\pla}{\hspace{1.5cm}}
\newcommand{\hz}{\vspace{0.5cm}}

\newcommand{\si}{\sigma}

\newcommand{\eps}{\varepsilon}

\newcommand{\noo}{\left \|}
\newcommand{\rrm}{\right \|}

\newcommand{\bet}{\left |}
\newcommand{\rag}{\right |}

\newcommand{\intt}{\int\limits}
\newcommand{\summ}{\sum\limits}

\renewcommand{\i}{\subset}

\documentstyle[12pt]{article}
\begin{document}

\title{ \bf Proportional subspaces of spaces with unconditional basis have
good
volume properties }

\author{ Marius Junge }

\date{}

\maketitle

\begin{abstract}
\hspace{-2.0em} A generalization of Lozanovskii's result is proved. Let E be
$k$-dimensional
subspace of an $n$-dimensional Banach space with unconditional basis.
Then there exist $x_1,..,x_k \subset E$ such that $B_E \p \subset \p
absconv\{x_1,..,x_k\}$
and
\[
\kla \frac{{\rm vol}(absconv\{x_1,..,x_k\})}{{\rm vol}(B_E)}
\mer^{\frac{1}{k}} \kl \kla
e\p \frac{n}{k} \mer^2 \pl .\]
This answers a question of V. Milman which appeared during a GAFA seminar talk
about the hyperplane
problem. We add logarithmical estimates concerning the hyperplane conjecture
for proportional subspaces and quotients
of  Banach spaces with unconditional basis.
\end{abstract}

\setcounter{lemma}{0}
\section*{Introduction}
An open problem in the theory of convex sets is the following \hz

{\bf Hyperplane problem:} {\it Does there exist a universal constant $c > 0$
such that for
all $\nen$ and all convex, symmetric bodies $K \i \rz^n$ one has}
\[ |K|^{\frac{n-1}{n}} \kl c\pl \sup_{H \p hyperplane} |K\cap H|  \pl ?\]

For some classes of convex sets there is a positive solution to this problem.
For example Bourgain first proved the existence of a constant independent of
dimension for the class of convex sets with unconditional basis. This can be
formulated as follows

\begin{theorem}[Bourgain] For all convex, symmetric bodies $K \i \rz^n$ one
has
\[ |K|^{\frac{n-1}{n}} \kl 2\pl  \sqrt{6} \pll \inf\{\kla \frac{|B|}{|K|}
\mer^{\frac{1}{n}}
\mitt K \i B \p \mbox{and } B \p\p \mbox{with unc. basis} \} \pll
\sup_{H \p hyperplane} |K\cap H|
\pl .\]
\end{theorem}

For further positive solutions and background information we refer to the
papers
of Ball \cite{BA}, Milman/Pajor \cite{MIPA} and the author \cite{JU}.
In a seminar talk about the hyperplane problem V. Milman asked
whether the unit ball of a proportional
subspaces of a Banach space with unconditional basis is well contained
(in the volume
sense) in a convex body with unconditional basis, more precisely, whether the
infimum on the right hand side of Bourgain's theorem is uniformly bounded for
proportional subspaces of Banach spaces with unconditional  basis. This can be
answered in the
positive.

\begin{theorem} Let X be a $n$-dimensional Banach space with unconditional
basis and E a $k$-dimensional subspace. Then there exist $x_1,..x_k \in E$
such that
\[ B_E \pl \subset \pl absconv\{x_1,..,x_k\} \pla  \mbox{and} \pla
\kla \frac{|absconv\{x_1,..,x_k\}|}{|B_E|} \mer^{\frac{1}{k}} \kl \kla
e\p \frac{n}{k} \mer^2 \pl .\]
\end{theorem}

This theorem is a generalization of Lozanovskii's result, which corresponds to
the case $k=n$. In fact we use his approach. In
particular, the above theorem gives a uniform bound for the hyperplane
problem in the case of proportional subspaces of a Banach space with
unconditional
basis. This includes proportional
subspaces of $\ell_{\infty}^n$ which are often used to produce more
or less pathological phenomena in the local theory of Banach spaces. For the
hyperplane problem the estimates of theorem 2 can even be improved to a
logarithmical order.

\begin{theorem} Let E be a $k$-dimensional subspace of a $n$-dimensional
Banach space with unconditional basis. Then one has

\[ |B_E|^{\frac{k-1}{k}} \kl 2e\pl \sqrt{6+3\ln\frac{n}{k}}
\pl  \sup_{H\p hyperplane} |B_E\cap H|
  \pl.\]

\end{theorem}

Apart from the geometric interpretation, a convex polytope with not to many
faces nearly satisfies the hyperplane conjecture, theorem 3 destroys the hope
of producing counter examples by taking `bad' subspaces of `good' spaces.
For convex polytopes with not too many extreme points, we can proof a
slightly  weaker result. Although in this case the operator ideal theory
which is involved in the proof is a little bit harder.


\begin{theorem} Let E be a $k$-dimensional quotient of an $n$-dimensional
Banach space with unconditional basis. Then one has

\[ |B_E|^{\frac{k-1}{k}} \kl c_0 \pl (1+\ln n) \pl
  \sup_{H \p hyperplane} |B_E\cap H|
\pl,\]
where $c_0$ is a universal constant.

\end{theorem}

\setcounter{lemma}{0}
\section*{Proofs}
We will use standard Banach space notation, in particular we denote by $B_X$
the
unit ball of a Banach space X. In contrast to this $B_p^n$ is the unit ball of
the
classical sequence space
$\ell_p^n$, $1 \le p \le \infty$. For the volume of a convex body $B \subset
\rz^n$
we use $\bet B \rag$. The same notation is used for the lower dimensional
volumes of sections
of a convex body. A Banach space X has a $(1-)$ unconditional basis if there
exists a basis
$(e_i)_{i\in I}$ such that for all signs $(\eps_i)_{i \in  I}$ and
coefficients $(\alpha_i)_{i \in  I}$
\[ \noo \summ_{i \in I} \eps_i \p \alpha_i \p e_i \rrm \pl \le \pl
   \noo \summ_{i \in I} \alpha_i \p e_i \rrm \pl .\]
The following lemma of Lozanovskii \cite{LO} is crucial for the following.

\begin{lemma} Let X be an n-dimensional Banach space with unconditional basis
$(e_i)_1^n$. Then there exists positive weights $(\lambda_i)_1^n$ such that
\[ \frac{1}{n} \summ_1^n \bet\alpha_i \rag \pl \le
\pl \noo \summ_1^n \alpha_i \lambda_i e_i \rrm \pl \le \pl \sup_{i=1;..,n}
\bet
\alpha_i\rag \pl .\]
\end{lemma}

The next lemma reduces the problem to subspaces of $\ell_1^n$.

\begin{lemma} Let X be a $n$-dimensional Banach space with unconditional
basis.
Then there exists an operator
$T : X \nach \ell_1^n$  with $\noo T \rrm \p \le \p1$ such that for every
$k$-dimensional subspace E one has
\[\kla \frac{\bet T^{-1}(B_1^n) \cap E \rag}{\bet B_E \rag} \mer^{\frac{1}{k}}
\pl \le \pl e \pl \frac{n}{k} \pll .\]
\end{lemma}

{\bf Proof:} Using the weights from lemma 1 we define
\[ T : X \nach \ell_1^n;\pl T(\summ_1^n \alpha_i\p e_i) \pl:=\pl
\kla\frac{\alpha_i}{n \lambda_i}\mer_1^n
\pll \mbox{and} \pll S : \l_{\infty}^n \nach X; \pl S( (\alpha_i)_1^n)\pl
:=\pl \summ_1^n n \lambda_i \p
 \alpha_i \p e_i \pl. \]
According to lemma 1 we have $\noo T \rrm \p \le \p 1$ and $\noo S \rrm \p \le
\p n$. For the subspace
$H \p :=\p T(E) \i \rz^n$ we can use Meyer/Pajor's volume estimate \cite{MEPA}
to deduce
\for
\kla \frac{\bet T^{-1}(B_1^n)\cap E \rag}{\bet B_E \rag} \mer^{\frac{1}{k}}
&=& \kla \frac{\bet H \cap B_1^n \rag}{\bet T(B_E) \rag} \mer^{\frac{1}{k}}
\pl = \pl
\kla \frac{\bet H \cap B_{\infty}^n\rag}{\bet T(B_E) \rag} \mer^{\frac{1}{k}}
\p \kla \frac{\bet H \cap B_1^n\rag}{\bet H \cap B_{\infty}^n \rag}
\mer^{\frac{1}{k}}\\
&\le& \kla \frac{\bet S(H\cap B_{\infty}^n) \rag}{\bet B_E \rag}
\mer^{\frac{1}{k}} \pl
\kla \frac{\bet B_1^k \rag}{\bet B_{\infty}^k \rag} \mer^{\frac{1}{k}}\\
&\le& n \pl (k!)^{-\frac{1}{k}} \pl \le \pl e \pl \frac{n}{k} \pll .\\[-1.2cm]
\mel \qed \hz

{\bf Proof of theorem 2:} By lemma 2 we are left to prove the assertion for a
$k$-dimensional subspaces H of $\ell_1^n$. For this let us denote by $P$ the
orthogonal
projection from $\ell_2^n$ onto H. Define $x_i \p:=\p P(f_j)$, where
$(f_i)_1^n$
denotes the standard unit vector basis in $\rz^n$. The polar of $H \cap B_1^n$
is a zonotope whose volume can be estimated with a well known determinant
formula
\cite{MCM}, namely
\for
|(B_1^n\cap H)^{\circ}|
&=& 2^k \pl \summ_{{\rm card}(\si)=k} \bet det_k(x_j)_{j\in \si}\rag \\
&\le& \kla {n \atop k} \mer \pl \sup_{{\rm card}(\si)=k} 2^k \bet
det_k(x_j)_{j \in \si}\rag \pl.
\mel
Now fix a subset $\si \subset \{1,..,n\}$ of cardinality $k$ where the
supremum
is attained (in particular the vectors $(x_j)_{j \in \si}$ are independent).
Clearly we have for all $x \in H$
\[ \noo x \rrm_1 \pl = \pl \summ_1^n \bet \langle x,x_j\rangle \rag \pl \ge
\pl
\summ_{j \in \si} \bet\langle x,x_j\rangle \rag \pl =: \pl \noo x \rrm_{\si}
\pl. \]
The unit ball $B_{\si}$ of the norm $\noo\p \rrm_{\si}$ is the image of an
$\ell_1^k$-ball
and contains $B_1^n \cap H$. By the inverse Santal\'{o} inequality for
zonoids,
due to Reisner \cite{RE}, we obtain
\for
|B_{\si}| \p |B_{\si}^{\circ}|
&=&  |B_1^k|\p |B_{\infty}^k| \pl \le \pl
  |B_1^n \cap H| \p |(B_1^n \cap H)^{\circ}|\\
&\le& |B_1^n \cap H|\p \kla {n \atop k}\mer \p |B_{\si}^{\circ}|\pl .
\mel
Therefore we have proved
\for \kla \frac{\bet B_{\si}\rag}{\bet B_1^n \cap H  \rag}\mer^{\frac{1}{k}}
&\le& \kla {n \atop k} \mer^{\frac{1}{k}} \pl \le \pl e \pl \frac{n}{k}
 \pll .\\[-1.0cm]
\mel \qed \hz

\begin{rem}
By duality we obtain that the unit ball $B$ of a $k$-dimensional quotient
of a $n$-dimensional Banach space with unconditional basis contains the
affine image of a cube $C$ with
\[ \kla \frac{|B|}{|C|} \mer^{\frac{1}{k}} \kl c_0 \pl \kla \frac{n}{k} \mer^2
\pl .\]
\end{rem}\hz

For the hyperplane problem let us recall that a symmetric, convex body $K$
is in isotropic position if
\begin{enumerate}
\item[i)] $|K|\pl =\pl 1$,
\item[ii)] $\intt_{K} \langle x,e_j\rangle \langle x,e_i\rangle dx \pl = \pl
L_K^2\p \delta_{ij}$.
\end{enumerate}
In this case $L_K$ is the constant of isotropy of $K$. Let us note that for
every convex, symmetric body there is an affine image which is in isotropic
position. With the help of this it's essentially
Hensley's result \cite{HEN}, that an upper bound for the constant
of isotropy solves the hyperplane problem for any position of $K$. For further
information see for instance \cite{MIPA}. In the following we will denote by
$E_K$ the
Banach space $\rz^n$ equipped with the gauge $\noo \p \rrm_K$, i.e. $E_K$
is the Banach space whose unit ball is $K$. It was already discovered
by K. Ball that the notion of (absolutely) p-summing ($1\kll p < \infty$)
is a useful tool
for certain estimates of the constant of isotropy. An operator $T :X \nach Y$
is p-summing if there exists a constant $c \gll 0$ such that for all $\nen$,
$(x_k)_1^n \subset X$
\[ \kla \summ_1^n \noo Tx_k \rrm^p \mer^{\frac{1}{p}} \kl c \pl
\sup_{\noo x^* \rrm_{X^*} \kll 1} \kla \summ_1^n \bet\langle x_k,x^* \rangle
\rag^p \mer^{\frac{1}{p}} \pl .\]
The best possible constant $c$ will be denoted by $\pi_p(T)$.

\begin{lemma} \label{p1} Let $K \i \rz^n$ be in isotropic position. For the
formal
identity $\iota: \ell_2^n \nach E_K$ one has
\[ L_K \p \pi_1(\iota^*) \kl 2\sqrt{2} \pl . \]
\end{lemma}

{\bf Proof:} As a consequence of C. Borell's lemma we have for all $\alpha\in
\rz^n$

\[ L_K \p\noo\alpha \rrm_2 \kl 2\sqrt{2} \pl \intt_K \bet \langle x,\alpha
\rangle \rag \p dx \pl .\]
(For the precise constant see \cite{MIPA}.) Now let $m \in \nz$,
$(\alpha_j)_1^m \subset \rz^n$. Then we have
\for
L_K\pl  \summ_1^m \noo \alpha_j \rrm_2 &\le& 2 \sqrt{2} \pl
 \summ_1^m \intt_K \bet\langle x,\alpha_j \rangle \rag \p dx\\
&=&  2\sqrt{2} \pl \intt_K \summ_1^m \bet\langle \frac{x}{\p\noo x
\rrm_K},\alpha \rangle \rag \p \noo x \rrm_K\p dx\\
&\le& 2 \sqrt{2} \pl  \intt_K \noo x \rrm_K \p dx \pl \sup_{\noo y \rrm_K \kll
1}
\summ_1^m  \bet\langle y,\alpha_j \rangle \rag \\
&\le& 2 \sqrt{2}\pl \sup_{\noo y \rrm_K \kll 1} \summ_1^m  \bet\langle
y,\alpha_j \rangle \rag \pl .\\[-1.3cm]
\mel \qed \hz

{\bf Proof of theorem 3:} Let E be a $k$-dimensional subspace of a
$n$-dimensional
Banach space X with unconditional basis. We can find an isotropic position
for the unit ball of E, i.e. there exists a linear map $T:\rz^k \nach X$ such
that
$E \p = T(\rz^k)$ and $K\p=\p T^{-1}(B_X)$ is in isotropic position. Let us
define $S \p :=\p T \iota :\ell_2^k \nach X$. By lemma \ref{p1} we have
\[ L_K \p \pi_1 (S^*) \kl 2 \sqrt{2} \pl .\]
Since $X$ has an unconditional basis the same is true for $X^*$ and therefore
$S^*$ well-factors through $\ell_1^n$ \cite[Lemma 8.15]{PS}. By duality there
exist $W :\ell_2^k \nach \ell_{\infty}^n$,
$\noo W \rrm \kll 1$ and
$V: \ell_{\infty}^n \nach X$ such that $S \p =\p VW$ and
\[ \noo V \rrm \p \kl 2 \sqrt{2} \pl L_K^{-1} \pl .\]
Let $B\p:=\p W^{-1}(B_{\infty}^n)$. From $Im(VW)\p=\p E$ we deduce
\[ S(B) \pl =\pl V(B_{\infty}^n \cap W(\ell_2^k)) \pl \subset \pl \noo V \rrm
\pl B_E \pl .\]
and therefore $B \subset \noo V \rrm K$. Gluskin's theorem together
with $\noo W \rrm \kll 1$ implies a lower estimate for the volume ob $B$.
Using the constant from \cite{BAPA} we obtain
\for
2 \sqrt{2} &=& 2\sqrt{2} \pl \kla \frac{|K|}{|B|} \mer^{\frac{1}{k}} \pl
     \kla \frac{|B|}{|K|} \mer^{\frac{1}{k}}\\
&\le& 2 \sqrt{2}\pl |B|^{-\frac{1}{k}} \pl \noo V \rrm \\
&\le& e \sqrt{2+\ln \frac{n}{k}} \pll 2\sqrt{2}\pl L_K^{-1} \pl .
\mel

This means $L_K \kll e \p \sqrt{2+\ln\frac{n}{k}}$. Hensley's theorem
\cite{HEN} yields the
assertion. \qed\hz

The logarithmic estimate of the hyperplane constant for quotient spaces is
based
on the use of C. Borell's lemma in a similar setting as in lemma \ref{p1}.

\begin{lemma}\label{p2} Let $K \i \rz^k$ be in isotropic position and $T:E_K^*
\nach Y$ an isometric embedding of $E_K^*$ in a $n$-dimensional Banach space
Y.
Then there exists an extension $S:Y\nach \ell_2^k$ of the formal
identity $\iota^* E_K^* \nach \ell_2^k$ with $ST \p=\p \iota_K^*$ and
\[ L_K \pl \pi_1 (S) \kl c_0 \pl \kla 1 + \ln n \mer \pl .\]
\end{lemma}

{\bf Proof:} Let $K$ be in isotropic position and denote by $\mu$ the Lebesgue
measure restricted on $K$. Choosing $p\p=\p2+\ln n \gll 2$
we want to construct a suitable factorization of $L_K\p \iota^*$.
For this consider $J: E_K^* \nach L_{\infty}(K,\mu)$, $ \alpha \mapsto \kla x
\mapsto \langle
x,\alpha\rangle \mer$. Clearly $\noo J \rrm \kll 1$. Since $L_{\infty}(K,\mu)$
has the
extension property, see \cite{PI1}, there is an operator $L : Y \nach
L_{\infty}(K,\mu)$ with $LT \p=\p J$ and $\noo L \rrm \kll 1$. Furthermore,
we define $I :L_{\infty}(K,\mu) \nach L_{p'}(K,\mu)$ the formal identity,
$p'$ the conjugate index to $p$, and
$P : L_{p'}(K,\mu) \nach \ell_2^k$ by
\[ P(f) \pl := \pl \kla \intt_{K}
f \frac{\langle x,e_j \rangle}{L_K} \p d\mu(x) \mer_1^k\pl.\]
It is easy to see
that $L_K \p \iota_K^* \p=\p PIJ$ and $S \p:=\p PIL$ is an appropriate
extension.
For the norm of P we deduce from C. Borell's lemma, see \cite[Appendix]{MIS},
and the isotropic position of $K$
\for
\noo P \rrm\!\! &=&\! \!\noo P^* \rrm\!\! \pl = \pl\!\!
 \sup_{\noo \beta \rrm_2 \kll 1} \kla \intt_K \bet \langle \frac{\beta}{L_K},x
 \rangle \rag^p \p dx \mer^{\frac{1}{p}}
\pl\!\! \le \pl\!\!  c_0 \pl\! p \!\pl \sup_{\noo \beta \rrm_2 \kll 1} \kla
\intt_K \bet \langle \frac{\beta}{L_K},x \rangle \rag^2 \p dx
\mer^{\frac{1}{2}}\!\! \kl\! \!c_0\! \pl p\pl\! .
\mel
In fact we have proved $\iota_{p'}(S) \kll c_0\p p$, where $\iota_{p'}$
denotes
the $p'$-integral norm. By the choice of $p$ the proof of the lemma will be
completed
if we can show \vspace{0.2cm}

{\bf $(*)$}\hfill $\pi_1(S) \kl n^{\frac{1}{p}}\pl \iota_{p'}(S)\pl .$
\hspace*{\fill}\vspace{0.2cm}

Given a sequence $(y_j)_1^m \subset Y$ with $\sup_{y^* \in B_{Y*}} \summ_1^m
\bet \langle y_j,y^* \rangle \rag \kll 1$ we define the operator $R:
\ell_{\infty}^m \nach Y$,
$R(\beta_j)_1^m \p = \p \summ_1^m \beta_j y_j$ whose norm is less than 1.
In this situation we can use an interpolation formula \cite{??}
for the $p$-summing norm to deduce
\[ \pi_p(R) \kl \pi_2(R)^{\frac{2}{p}} \p \noo R \rrm^{1-\frac{2}{p}}
\kl n^{\frac{1}{p}}\pl. \]
Here we have used the well-known fact $\pi_2(R) \kll \sqrt{n} \p \noo R \rrm$
for any operator of rank at most $n$, see for example \cite{PI2}.
Now we can find $(\alpha_j)_1^m
\subset B_2^k$ with $\noo S(y_j) \rrm \p=\p \langle S(y_j), \alpha_j
\rangle$.
Clearly the operator $V :\ell_2^k \nach \ell_{\infty}^m$, $V(x)\p:=\p (\langle
x,\alpha_j \rangle)_1^m$
has also of norm at most $1$ and trace duality (see for example \cite{PI1})
implies
\for
\summ_1^m \noo S(y_j) \rrm &=& tr(VSR) \kl \iota_{p'}(VS) \pl \pi_p(R) \kl
\iota_{p'}(S) \pl n^{\frac{1}{p}}\pl .
\mel
By the definition of the $\pi_1$-summing norm we have proved $(*)$.
\hfill $\Box$\hz

\begin{rem}
In the proof above an isometric embedding is not really needed. The
$\pi_1$-summing norm of an extension can be chosen according to the minimal
distance
of $E_K^*$ to a $k$-dimensional subspace of $Y$.
\end{rem}

Given lemma \ref{p2} the proof of theorem 4 of the introduction follows the
same pattern
as the proof of theorem 3.\hz

{\bf Proof of theorem 4:} Let $X$ be a $n$ dimensional Banach space with
unconditional basis.
For a $k$-dimensional quotient space $E$ of $X$ with quotient map $Q: X \nach
E$
we can find an isomorphism $I :E \nach \rz^k$, such that $K\p=\p I(B_E)$ is in
isotropic position. In this case $T := Q^* I^* : E_K^* \nach X^*$ defines
an isometric embedding. Applying lemma \ref{p2} there is an extension $S: X^*
\nach \ell_2^k$
of $L_K \iota^*$ with $\pi_1(S) \kll c_0 \p (1+\ln n)$. Since $X^*$ also an
unconditional basis
$S$ factors through $\ell_1^n$ \cite[Lemma 8.15]{PS}. More precisely, there
are $W: X^* \nach \ell_1^n$,
$\noo W \rrm \kll 1$ and $V: \ell_1^n \nach \ell_2^k$ with $S \p=\p VW$ and
\[ \noo V \rrm \kl \pi_1(S) \kl c_0 \pl (1+\ln n) \pl .\]
Now we consider the $k$-dimensional subspace $F \p:=\p  WT(E^*) \i \ell_1^n$.
Instead of Gluskin's estimate we can use a dual volume estimate first
essentially proved by Figiel and Johnson \cite{FIJ}
\[ \sqrt{k} \pll \kla \frac{|V(B_F)|}{|B_2^k|} \mer^{\frac{1}{k}} \kl
c_1 \noo V \rrm \pl .\]
(Indeed, $\ell_1$ is of cotype 2 and therefore every subspace has bounded
volume ratio. The inequality
follows from this if we note that by Grothendieck's theorem V is 2-summing.)
Since $ST\iota^* \p=\p L_K id_{\rz^k}$ we conclude with the inverse
of Santal\'{o}'s inequality \cite{BM}
\for
L_K&\le& L_K \pl \sqrt{k} \pl  \kla \frac{|B_2^k|}{|K|} \mer^{\frac{1}{k}}
   \kl c_2 \pl \sqrt{k} \pl  \kla \frac{|\stackrel{\circ}{K}|}{|B_2^k|}
   \mer^{\frac{1}{k}}\\
&=& c_2 \pl \kla \frac{|W(B_{E^*})|}{|B_F|} \mer^{\frac{1}{k}}\pl
\sqrt{k} \pl \kla \frac{|V(B_F)|}{|B_2^k|} \mer^{\frac{1}{k}}\\[+0.3cm]
&\le& c_2 \pl \noo W \rrm \pl c_1 \pl \noo V \rrm \kl c_0 \pl c_1 \pl c_2 \pl
(1+ \ln n)\pl .\\[+0.1cm]
\mel
Hensley's theorem implies the assertion,
see \cite{HEN} and \cite{JU}.\hfill $\Box$\hz

\begin{frem} For the proofs of theorem 3, 4  we have used operator ideal
techniques.
This allows us to formulate the results in a little bit stronger form which
is similar to the formulation of Bourgain's theorem in the introduction.
Let $K \subset \rz^k$ then we have \vspace{0.3cm}

\for
|K|^{\frac{k-1}{k}} &\le&  \inf\{\kla \frac{|B_E|}{|K|} \mer^{\frac{1}{k}}
\mitt E \subset X \mbox{with unc. basis and } dimX=n\}\\[+0.3cm]
& & \times \pl 2e \pl  \sqrt{6+3\ln\frac{n}{k}} \pl \sup_{H \p hyperplane}
|K\cap H|
\mel

and

\for
|K|^{\frac{k-1}{k}} &\le& \inf\{\kla \frac{|B_E|}{|K|} \mer^{\frac{1}{k}}
\mitt E \mbox{ quotient of } X \mbox{with unc. basis and } dimX=n \}\\[+0.3cm]
& & \times \pl c_0 \pl (1+\ln n) \pl \sup_{H \p hyperplane} |K\cap H| \pl.
\mel

\end{frem}


\hz
1991 Mathematics Subject Classification: 52A38, 46B45, 52A21.

Key words: Unconditional basis, volume, hyperplane conjecture.
\begin{quote}
Marius Junge

Mathematisches Seminar der Universit\H{a}t Kiel

Ludewig-Meyn-Str. 4

24098 Kiel

Germany

Email: nms06@rz.uni-kiel.d400.de
\end{quote}
\end{document}